\documentclass[12pt]{article}
\usepackage{MnSymbol}
\usepackage{amsmath,amsfonts,
verbatim}

\newcommand{\C}{{\mathbb C}}
\newcommand{\N}{{\mathbb N}}
\newcommand{\Z}{{\mathbb Z}}

\newcommand{\LL}{\mathcal{L}}
\newcommand{\Deck}{\mathrm{GDeck}}
\newcommand{\F}{\mathbb{F}}
\newcommand{\U}{{\mathbb U}}

\newcommand{\X}{{\mathbb X}}
\newcommand{\XX}{\tilde{\mathbf{X}}}
\newcommand{\Y}{\mathbb{Y}}

\newcommand{\Gal}{\mathrm{Gal}}

\newcommand{\pr}{{\rm pr}}

\newcommand{\MM}{\mathcal{M}}
 
\newcommand{\la}{\langle}
\newcommand{\ra}{\rangle}
\newtheorem{pkt}{}[section]  
\newcommand{\bpk}{\begin{pkt}\rm }  
\newcommand{\epk}{\end{pkt}} 
\newcommand{\inv}{^{-1}}   

\newcommand{\be}{\begin{equation}}  
\newcommand{\ee}{\end{equation}}

\newcommand{\dcl}{\mathrm{dcl}}

\newcommand{\kk}{\mathrm{k}}

\newcommand{\pp}{\mathbf{p}}
\newcommand{\Aut}{\mathrm{Aut}}

\usepackage{graphicx}

\usepackage{tikz}
\usetikzlibrary{arrows}
\tikzstyle{block}=[draw opacity=0.7,line width=1.4cm]

\title{Section and towers}
\author{B.Zilber\footnote{Supported by the EPSRC program grant ``Symmetries and Correspondences''}}

\begin{document}
\maketitle

\abstract{We discuss  the towers of finite \'etale covers which were essentially introduced by A.Tamagawa \cite{Ta} and used e.g. in \cite{Sa}. The statement about correspondence between sections and cofinal towers is a folklore but perhaps not in a very explicit form. The last section explains how the "injectivity statement" of Grothendieck section conjecture fails for abelian varieties, which is also known in some form from \cite{CStix}.

The paper is based on \cite{AZ} which was aimed to reinterpret anabelian setting in model theory terms.}

 \section{A short overview of structure $\XX^{et}$}\label{ss2}
We start with an overview of the key structure introduced and studied in \cite{AZ}. It is essentially the projective object - the Grothendieck universal \'etale cover of a smooth $\kk$-variety $\X.$

\begin{center}
 \begin{tikzpicture}
 \draw (9.8,8.5) node {$\tilde{\X}$};  
 \draw[ ->][line width=0.3mm] (9.6,  8.7) arc (60:300:0.2);
\draw (9.1,8.5) node {$  {\Gamma}$};
 \draw[ ->][line width=0.5mm] (9.8,  8.2)--(10,0.4); 
\draw (10.1,4.1) node {$ \pp$};

\draw (12.3,3.4) node {$\tilde{ \pp}_{\mu \alpha_n}$};
\draw (11.3,3.4) node {$\tilde{ \pp}_{\mu \alpha_1}$};

\draw (11,2) node {$\X_{\mu\alpha_1}$};
 \draw (12.5,2) node {$\X_{\mu\alpha_n}$};
 \draw (11.6,2) node {$\ldots$};

\draw[->][line width=0.3mm] (9.8,8 )--(11, 2.5);
\draw[->][line width=0.3mm] (9.8,8 )--(10.8, 2.5);
\draw[->][line width=0.3mm] (9.8,8 )--(10.9, 2.5);
 \draw (11.5,2.6) node {$\ldots$};
\draw[->][line width=0.3mm] (9.8,8 )--(12.0, 2.5);
\draw[->][line width=0.3mm] (9.8,8 )--(12.2, 2.5);
\draw[->][line width=0.3mm] (9.8,8 )--(12.1, 2.5);

\draw (8,5) node {$\tilde{\pp}_{\nu\beta_1}$};
 \draw (9.4,5) node {$\tilde{\pp}_{\nu\beta_k}$};
 
\draw (8.2,3.5) node {$\X_{\nu\beta_1}$};
\draw (9.6,3.5) node {$\X_{\nu\beta_k}$};
\draw (8.9,3.5) node {$\ldots$};
\draw[ ->][line width=0.3mm] (7.6, 3.7) arc (60:300:0.2);
\draw (6.9,3) node {$\mathrm{GDeck}(\X_{\nu\beta_1}/\X)$};

\draw[->][line width=0.3mm] (9.8,8 )--(8.2, 4);
\draw[->][line width=0.3mm] (9.8,8 )--(8.3, 4);
\draw[->][line width=0.3mm] (9.8,8 )--(8.4, 4);
 \draw (9,4.3) node {$\ldots$};
\draw[->][line width=0.3mm] (9.8,8 )--(9.5, 4);
\draw[->][line width=0.3mm] (9.8,8 )--(9.7, 4);
\draw[->][line width=0.3mm] (9.8,8 )--(9.6, 4);



\draw[->][line width=0.5mm](10.8,1.8)--(10.1,0.4);
\draw[->][line width=0.5mm](12.3,1.8)--(10.2,0.4);
\draw (11,1.4) node {$\mu_{\alpha_1}$};
\draw (12.2,1.4) node {$\mu_{\alpha_k}$};

\draw[->][line width=0.5mm](9.5,3.2)--(10,0.4);

\draw[->][line width=0.5mm](8.3,3.2)--(9.9,0.4);
\draw (8.8,1.8) node {$\nu_{\beta_1}$};
\draw (9.8,1.8) node {$\nu_{\beta_k}$};

\draw (10,0) node {$\X$};

\draw[->][line width=0.1mm] (9.6,3.2 )--(10.8, 2.3);
\draw[->][line width=0.1mm] (9.6,3.2 )--(10.7, 2.3);

\draw[ ->][line width=0.3mm] (12.8, 1.8) arc (220:480:0.2);
\draw (14.7,1.8) node {$\mathrm{GDeck}(\X_{\mu\alpha_n}/\X)$};


\end{tikzpicture}

The diagram for $\XX^{et}.$

\bigskip

\end{center}
{\bf Explaining the picture} (see \cite{AZ}, 7.1-7.3 and Corollary 7.11)
\bpk All arrow diagrams commute.
\epk
\bpk
Each $\X_{\nu,\beta_i}$ is an absolutely irreducible variety over the field $\kk[\beta_1]=...=\kk[\beta_k],$  a Galois extension of $\kk.$
 $ \X_{\nu,\beta_i}(\kk^{alg})$ is the set of its $\kk^{alg}$-points, a subset of a projective space.
\epk
\bpk
\label{item[3]} Each $\nu_{\beta_i}$  is an \'etale covering map $\nu_{\beta_i}: \X_{\nu,\beta_i}(\kk^{alg})\to \X(\kk^{alg})$. 
\epk
\bpk
\label{item[4]} $ \tilde{\X}(\kk^{alg})$ is a set with the regular action of a  group
$   {\Gamma}.$
\epk
\bpk
\label{item[5]}
Each $\tilde{\pp}_{\nu\beta_i}$ is a finite collection of surjective {\bf maps}
$$p: \tilde{\X}(\kk^{alg})\to \X_{\nu\beta_i}(\kk^{alg}).$$
In particular, if  $\X_{\nu\beta_i}(\kk^{alg})= \X_{\mu\alpha_j}(\kk^{alg})$ and $\nu_{\beta_i}=\mu_{\alpha_j}$
then  $\tilde{\pp}_{\nu\beta_i}=\tilde{\pp}_{\mu\alpha_j}.$\\
 In case $ \X_{\nu,\beta_i}(\kk^{alg})= \X(\kk^{alg})$ the collection
$\tilde{\pp}_{\nu\beta_i}$ consists of one map $\pp.$

\epk
\bpk \label{item[5+]} Suppose there is a morphism $(\mu_\alpha\inv\nu_\beta):\X_{\nu,\beta}(\kk^{alg})\to \X_{\mu,\alpha}(\kk^{alg})$ of \'etale covers (see notation in \cite{AZ}, section 4) . Then for every  $p\in \tilde{\pp}_{\mu,\alpha}$ there is $q\in \tilde{\pp}_{\nu,\beta}$ such that \be \label{after13}
(\mu_\alpha\inv\nu_\beta)=p\circ q\inv,\ee
and for every $q\in \tilde{\pp}_{\nu,\beta}$ there is $p\in \tilde{\pp}_{\mu,\alpha}$ such that (\ref{after13}) holds.

\epk
\bpk
\label{item[6]}
Given $p\in\tilde{\pp}_{\nu\beta_i}$ 
$$\tilde{\pp}_{\nu\beta_i}=\{ g\circ p: g\in \mathrm{GDeck}(\X_{\nu\beta_i}/\X)\}$$
 where $\mathrm{GDeck}(\X_{\nu\beta_i}/\X)\}$ is the geometric deck-transformation group. 
\epk
\bpk
\label{item[7]} The fibres of  $\pp$ are
$   {\Gamma}$-orbits. The fibres of  $p\in \tilde{\pp}_{\nu,\beta_i}$ are orbits by a finite index normal subgroup 
$  {\Delta}_{\nu,\beta_i}$ of $   {\Gamma}.$ 
$$\mathrm{GDeck}(\X_{\nu\beta_i}/\X)\cong    {\Gamma}/  {\Delta}_{\nu,\beta_i}.$$
\epk
\bpk
\label{item[8]}  For each finite collection  $\X_{\lambda_1\gamma_1}, \ldots\X_{\lambda_m\gamma_m},$ there is a $\nu_\beta$ such that
$$  {\Delta}_{\nu,\beta}\le \bigcap_{0<j \le m}  {\Delta}_{\lambda_j,\gamma_j}.$$ 
\epk
\bpk
 \label{item[9]} $$  \bigcap_{\mbox{ all } \nu,\beta} \Delta_{\nu,\beta}= \{ 1\}.$$ 
\epk
\bpk
\label{item[10]}
$$\Aut\, \XX^{et}(\kk^{alg})\cong \pi_1^{et}(\X,x)$$
 
 \epk
\section{Sections and towers}
\bpk Let $$T(\X):\ \ 
\X \leftarrow \X_1\leftarrow \X_2 \leftarrow\ldots \X_i\leftarrow \X_{i+1}\leftarrow \ldots $$ be a tower of smooth complex algebraic varieties and
unramified  covers, all defined over $\kk.$ Let $$\Gamma_{T(\X)}=\lim_\leftarrow\Deck(\X_i/\X).$$
We call the tower {\bf cofinal} if $$\Gamma_{T(\X)}\cong \hat{\pi}_1^{top}(\X)$$
as profinite groups.
\epk 
\bpk \label{towers1}{\bf Proposition.} {\em Given $\X$  there is a cofinal chain $$\Gamma>\Delta_1>\ldots \Delta_i> \Delta_{i+1}>\ldots$$
of $\Aut \XX^{et}(\F)$-invariant  normal  finite index subgroups of  $\hat{\pi}_1^{top}(\X)=\Gamma.$

Given a section $s$ and a cofinal chain $\{ \Delta_i\}$ of $\Aut \XX^{et}(\F)$-invariant  normal  finite index subgroups of $\Gamma$ there exists a 
 tower $T_s(\X)$ over $\kk$ such that $$\Deck(\X_{i}/\X)\cong \Gamma/\Delta_i.$$
 }

{\bf Proof.} Let $\Gamma:=\hat{\pi}_1^{top}(\X).$ $s\Gal_\kk$ acts on 
 $\Gamma$ since group $\Gamma$ is definable in $\XX^{et},$
 In particular, $s\Gal_\kk$ acts on the set of all finite index subgroups.
 \medskip
 
 {\bf Claim 1.} 
  There exists a decreasing sequence $\{ \Delta_n: n\in \N\}$ (depending on $\X$ only)
of   
$\Aut \XX^{et}(\F)$-invariant 
normal subgroups of $\Gamma$ of finite index with $\bigcap_n \Delta_n=\{ 1\}.$

Proof.
For each $\mu\in \MM_\X$ consider the subgroup $\Delta_\mu< \Gamma$ 
$$\Delta_\mu=\{ \gamma\in \Gamma: \forall p\in \tilde{\pp}_\mu \forall u\in \U \  p^\gamma(u)=p(u)\}$$
where $p^\gamma$ is the map $u\mapsto p(\gamma\cdot u).$ 

By 4.15 of \cite{AZ} 
$$\Delta_\mu= \bigcup_{\alpha \in \mathrm{Zeros}\mathbf{f}_\mu}\Delta_{\mu,\alpha},$$ the intersection of subgroups of periods of the maps $p: \U\to \X_{\mu,\alpha}$ which are finite index.
Hence $\Delta_\mu$ is of finite index in $\Gamma.$
It also follows that the intersection of the $\Delta_\mu$ is trivial.
It remains to choose a linearly ordered  cofinal subset ${\Delta_n: n\in \N}$ in ${\Delta_\mu: \mu\in \MM_\X}.$ Claim proved.

Let \be \label{Un}
\U_n:= \Delta_n\backslash \U, \ \bar{\pp}_n:\U\to \X,\ \bar{\pp}_{n,m}: \U_n\to \U_m\ee
where $\bar{\pp}_n$ is the covering map induced by $\pp:\U\to \X$ on $\U_n$ (recall that fibres of $\pp$ are $\Gamma$-orbits) and 
$\bar{\pp}_{n,m}$ is the map induced by the embedding $\Delta_n\le \Delta_m.$ 

Note that the $\U_n,$ $\bar{\pp}_n$ and $\bar{\pp}_{n,m}$ are $\Aut \XX^{et}(\F)$-invariant and so the action of $s\Gal_\kk$ on $\U$ induces the action on the tower 
$$\U_1\leftarrow\U_2\leftarrow\ldots$$

{\bf Claim 2.} The $\U_n$ can be given structure of smooth projective algebraic varieties defined over $\kk.$

Proof. By the argument in the proof of Claim 1, $\Delta_n=\Delta_{\mu,\alpha}=\mathrm{Per}\, p,$ for some $\mu_\alpha,$ $p: \U\to \X_{\mu,\alpha}.$ Set $p_n: \U_n\to  \X_{\mu,\alpha}$ be the bijective map induced by $p$ on $\U_n.$
We may assume that the set    $ \X_{\mu,\alpha}(\F)$ and the map $p$ are $s\Gal_{\kk[\alpha]}$-invariant, by possibly extending $\kk[\alpha]$ without changing the set and the map. Call $i_{n,\alpha}$ the map $p_n\inv:  \X_{\mu,\alpha}(\F)\to \U(\F).$ Note that by applying Galois conjugation we obtain a finite family $$\{ i_{n,\alpha}:  \X_{\mu,\alpha}(\F)\to \U(\F); \ \alpha\in \mathrm{Zeros} \mathbf{f}_\mu\}$$
of bijections.

Let $$\Y_n=\{ \la x,\alpha\ra: x\in  \X_{\mu,\alpha}\ \& \  \alpha\in \mathrm{Zeros} \mathbf{f}_\mu\}$$
the disjoint union of  $\kk[\alpha]$-varieties isomorphic to   $\X_{\mu,\alpha}.$ Let $i_n: \Y\to \U_n$ be the surjective map defined
as $$i_n(y)=u \leftrightarrow \exists \alpha \exists x\in  \X_{\mu,\alpha}\
y=\la x,\alpha\ra\ \&\
i_{n,\alpha}(x)=u.$$
By construction $\Y_n$ and $i_n$ are $\Gal_\kk$-invariant. 

Let $G$ be the group $\Gal( \kk[\alpha]:\kk)$ (recall that by our assumptions $\kk[\alpha]:\kk$ is Galois. For each $u\in \U_n,$
define the action of $G$ on $i_n\inv(u).$ Note that by construction
$$i_n\inv(u)=\{ \la x_\alpha, \alpha\ra:\ \alpha\in \mathrm{Zeros} \mathbf{f}_\mu\}$$
for some $x_\alpha\in  \X_{\mu,\alpha}.$ For $\sigma\in G$ set
$$\sigma: \la x_\alpha, \alpha\ra\mapsto \la x_{\sigma(\alpha)}, \sigma(\alpha)\ra.$$

By construction $G\backslash \Y_n$ is in bijective $\kk$-definable correspondence with $i_n(\Y_n)$ that is with $\U_n,$ that is
$$\U_n\cong   G\backslash \Y_n.$$
The object on the right is the quotient of smooth projective variety (reducible, in general)
by a regular action of a finite group. Hence $G\backslash \Y_n$   is isomorphic\footnote{Reference?}
 to a smooth projective variety $\X_n$ over $\kk$ via a surjective map $t_n: \Y_n\to \X_n$ with fibres which are $G$-orbits.
 Thus
there is a $s\Gal_\kk$-invariant bijective map onto the $\kk$-variety 
$$\mathbf{i}_n: \U_n\to \X_n.$$
Claim proved.

\medskip

Note that $t_{n,\alpha},$ the restriction of $t_n$ to $\X_{\mu,\alpha}\times \{ \alpha\},$ a component of $\Y_n,$ is a biregular isomorphism $t_{n,\alpha}:\la x,\alpha\ra \mapsto G\cdot  \la x,\alpha\ra$ on $\X_n$ defined over $\kk[\alpha].$ Consider the map
 $$t'_{n,\alpha}: x \mapsto G\cdot  \la x,\alpha\ra,\ \ \X_{\mu,\alpha}\to \X_n$$
which for simplicity of notation we call  $t_{n,\alpha}$ as well.
By construction
$$ t_{n,\alpha}\circ p_n= \mathbf{i}_n.$$
Define $\mathbf{j}_{n,m}: \X_n\to \X_m$ to be $\mathbf{j}_{n,m}=\mathbf{i}_m\circ\bar{\pp}_{nm}\circ \mathbf{i}_n\inv.$ This is definable over $\kk$ since  $\mathbf{i}_m,\ \bar{\pp}_{nm}$ and $ \mathbf{i}_n$ are $s\Gal_\kk$-invariant.
This is also a Zariski regular map
since by above
$$\mathbf{j}_{n,m}=t_{m,\beta}\circ p_m\circ \bar{\pp}_{n,m}\circ p_n\inv\circ t_{n,\alpha}\inv=t_{m,\beta}\circ (\nu_\beta\inv\mu_\alpha)\circ t_{n,\alpha}\inv  $$
where   $(\nu_\beta\inv\mu_\alpha): \X_{\mu,\alpha}\to \X_{\nu,\beta}$ is an intermediate regular map which can be presented as  $p_m\circ \bar{\pp}_{n,m}\circ p_n\inv.$

This gives us the cofinal tower 
$$T_s(\X): \ \X \leftarrow \X_1\leftarrow \X_2 \leftarrow\ldots \X_i\leftarrow \X_{i+1}\leftarrow \ldots $$
where the arrows $\X_{i+1}\to \X_i$ stand for the regular maps $\mathbf{j}_{i+1,i}.$ $\Box$
\epk
\bpk\label{Corol} {\bf Corollary} (of the proof). {\em 
Given $s$ and the tower $\{ \Delta_i: i\in \N\}$ of (\ref{Un}) the tower $T_s(\X)$ is determined uniquely up to isomorphism over $\kk.$ 

The system of bijections $\mathbf{i}_i$
$$\X \leftarrow \U_1\leftarrow \U_2 \leftarrow\ldots \U_i\leftarrow \U_{i+1}\leftarrow \ldots $$
$$ \downarrow \mathbf{i} \   \downarrow \mathbf{i}_1 \ \  \downarrow \mathbf{i}_2 \ldots\  \downarrow \mathbf{i}_i \ldots \downarrow \mathbf{i}_{i+1}\ldots $$ 
$$\X \leftarrow_{j_1} \X_1\leftarrow_{j_2} \X_2 \leftarrow\ldots \X_i\leftarrow_{j_{i+1}} \X_{i+1}\leftarrow_{j_{i+2}} \ldots $$
furnishes isomorphism between   the   structure on the
 tower of 
the $\U_i$ induced by the action of $s\Gal_\kk$
and the tower
$T_s(\X).$ 

Given any other such $s\Gal_\kk$-invariant tower
 $$T'_s(\X): \ \X \leftarrow \X'_1\leftarrow \X'_2 \leftarrow\ldots \X'_i\leftarrow \X'_{i+1}\leftarrow \ldots $$
  with covering maps $\mathbf{j}'_{i+1}:\X_{i+1}\to \X_i$
   there are isomorphism $q_i: \X_i\to \X'_i$ over $\kk$ such that
  $$q_i\circ \mathbf{j}_{i+1}=\mathbf{j}'_{i+1}\circ q_{i+1}.$$
  }
\epk
\bpk {\bf Proposition.} {\em Let $$\mathcal{T}(\X)=\{ T(\X): \{\Delta_i: i\in \N\}- \mbox{towers}\}$$
the set of  all  $ \{ \Delta_i: i\in \N\}$- towers  over $\kk.$\footnote{
That is $\Deck(\X_i/\X)\cong \Gamma/\Delta_i.$ 
Have to assume here that the tower
$\Deck(\X_l/\X)$ has unique, up to isomorphism of $\Gamma,$ presentation in the form $\Gamma/\Delta_i.$}   
Let $$\mathcal{S}(\X)=\{ s: \Gal_\kk\to \Aut\,\XX^{et}(\kk^{alg})\}$$
the set of all sections of $\pr: \Aut\,\XX^{et}(\kk^{alg})\to \Gal_\kk.$

Then the map $$s\mapsto  T_s(\X)$$
induces a bijection $$\mathcal{S}(\X)_{\slash\mbox{conj}} \to \mathcal{T}(\X)_{\slash \mbox{iso}}$$
between the set of  section modulo conjugation and the set of towers modulo isomorphisms over $\kk.$
}

{\bf Proof.} The map $s\mapsto  T_s(\X)_{\slash \mbox{iso}}$ is constructed above, see \ref{Corol}.
We construct the inverse map 
$$T(\X)_{\slash \mbox{iso}}\mapsto s_{\slash\mbox{conj}};\ \ \mathcal{T}(\X)_{\slash\mbox{iso}} \to \mathcal{S}(\X)_{\slash\mbox{conj}}.$$
Let $T(\X)$ be a $\Gal_\kk$-invariant  $ \{ \Delta_i: i\in \N\}$- tower. By the construction of $\XX^{et}$ the tower can be embedded into $\XX^{et},$ that is $\X_i=\X_{\mu_i,\alpha_i}$ for some $\mu_i\in \MM_\X,$ $\alpha_i\in \mathbf{f}_{\mu_i}$  and the $\mathbf{j}_{i+1}$ are
appropriate intermediate morphisms. Since the tower is over $\kk$ we can drop $\alpha_i.$ We also write $i$ for $\mu_i.$ 

Now we consider the respective sets of covering maps $p: \U\to \X_i, \ \ p\in \tilde{\pp}_i$ for each $i\in \N.$ 
 
 {\bf Claim 1.} There is a sequence $\pp_i\in \tilde{\pp}_i, i\in \N$ of covering maps
 $\pp_i: \U\to \X_i$ such that \be \label{pi}\mathbf{j}_{i+1}\circ  \pp_{i+1}=\pp_i, \mbox{ all }i\in \N.\ee
 
 Proof. By induction. For $i=0,$ set $\X_0:=\X$ and
 $\pp_0:=\pp.$ Suppose $\pp_n,$ $n\le i$ have been constructed satisfying   the requirement. We can choose $\pp_{i+1}$ by property \ref{item[5+]}. Claim proved
 
 {\bf Claim 2.} Suppose $\{ \pp'_i: i\in \N\}$ is another sequence satisfying (\ref{pi}). Then there is $\gamma\in \Gamma$ such that
 $$\pp'_i=\pp^\gamma_i, \mbox{ all }i\in \N,$$
 where $\pp^\gamma_i(u):=\pp_i(\gamma\cdot u)$ for all $u\in \U.$
 
 Proof. Choose $u\in \U$ and set $u_i:=\pp_i(u).$ First we prove that for each $n\in \N$
 there exists $u'\in \U$ such that $\pp'_i(u')=u_i$ for all $i\le n.$ 
  And for that it is enough to find $u'$ such that $\pp'_n(u')=u_n,$ since then 
 $$\pp'_{n-1}(u')=\mathbf{j}_{n}(\pp'_n(u'))=u_{n-1}, \ldots \pp'_{n-2}(u')=\ldots$$
  Note that $u'=u$ when $\pp'_0=\pp=\pp_0.$
  
By induction we assume that $\pp'_n(u')=u_n$
and need to find $u''$ such that $\pp'_{n+1}(u'')=u_{n+1}.$ Note that by (\ref{pi})
$\mathbf{j}_{n+1}(\pp'_{n+1}(u'))=u_n$ and so
$$\pp'_{n+1}(u')=g\cdot u_{n+1}\mbox{ for some }g\in \Deck(\X_{n+1}/\X_n).$$ 
We can find $\gamma\in \Gamma$ such that
$$\pp'_{n+1}(\gamma\inv u')=g\inv\pp'_{n+1}(u')=u_{n+1}.$$
Hence $u''=\gamma\inv u'$  satisfies the required.

Since the structure $\XX^{et}$ is compact in the profinite topology
there is an $u'$ which satisfies $\pp'_i(u')=u_i$ for all $i\in \N.$
Clealrly, $u$ and $u'$ are in the same fibre of $\pp$ and thus
$u'=\gamma\cdot u$ for some $\gamma\in \Gamma.$ Hence
$$\pp'_i(u)=\pp^\gamma_i(u)=g_i\cdot \pp_i(u).$$
It follows\footnote{Use the fact that groups of periods of both $\pp_i$ and $\pp'_i$ are $\Delta_i.$} that the equality holds for all $u.$ Claim proved.

{\bf Claim 3.} Any two sequences $\{ \pp_i\}$ and $\{\pp'_i\}$ satisfying (\ref{pi}) satisfy the same type over the sort $\F.$

Proof. By Claim 2 the sequence are conjugated by an element of $\gamma\in\Gamma.$ By construction the map $u\mapsto \gamma\cdot u$ is an automorphism of $\Aut\, \XX^{et}(\F)$ fixing all elements of sort $\F.$

{\bf Claim 4.} Let $\XX^{et}_ {  \{ \pp_i\}}(\F)$ be the structure
   $\XX^{et}(\F)$ with $\{ \pp_i\}$ named. The definable relation 
   on the sort $\F$ in the structure are exactly those which are definable in $\F_{| \kk},$ the field with constants for elements of $\kk.$
   
   Proof. By \cite{AZ}, Theorem 7.5,
   it is enough to prove that the definable relations on $\F$ in $\XX^{et}_ {  \{ \pp_i\}}(\F)$ are the same as in $\XX^{et}(\F).$ 
   
   Let $\varphi(\bar{x},  \pp_1,\ldots,\pp_n)$ be a formula in the language $\LL_\X(\{ \pp_i\})$ (the language of structure    $\XX^{et}_ {  \{ \pp_i\}}(\F)$), $\bar{x}$ a tuple of variables of sort $\F.$ 
 By Claim 3 there is a  formula $\psi_n(p_1,\ldots,p_n)$ in language $\LL_\X$ which is equivalent to a complete type of $\la \pp_1,\ldots,\pp_n\ra$ over $\F.$ We may assume that
 $$\varphi(\bar{x},  p_1,\ldots,p_n) \to \psi_n(p_1,\ldots,p_n).$$
 Now it is easy to see that in $\XX^{et}_ {  \{ \pp_i\}}(\F)$
 $$\varphi(\bar{x},  \pp_1,\ldots,\pp_n)\equiv
 \exists\   p_1,\ldots,p_n\ \psi_n(p_1,\ldots,p_n)\ \& \ \varphi(\bar{x},  p_1,\ldots,p_n)$$
 The formula on the right of $\equiv$ is in the language $\LL_\X$ and defines the relation $\varphi(\bar{x},  \pp_1,\ldots,\pp_n)$ in terms of $\XX^{et}(\F).$ Claim proved.
 
 {\bf Claim 5.} Let $T(\X)\in \mathcal{T}(\X)$ and $\{ \pp_i\}$ an associated sequence satisfying (\ref{pi}).
  Any automorphism $\sigma$ of $\F_{|\kk}$ induces a unique automorphism $s(\sigma)$ of 
$\XX^{et}_ {  \{ \pp_i\}}(\F).$ 

Proof. First note that 
$\sigma,$ being an automorphism of the field $\F,$ defines a transformation $\hat{\sigma}$  on algebraic sorts of $\XX^{et}_ {  \{ \pp_i\}}(\F),$ $$\hat{\sigma}: \X_{\mu,\alpha}(\F)\to \X_{\mu,\sigma(\alpha)}.$$
This transformation is an elementary monomorphism of $\XX^{et}_ {  \{ \pp_i\}}(\F),$ i.e. it
preserves the relation induced on the algebraic sorts in the structure $\XX^{et}_ {  \{ \pp_i\}}(\F).$ Indeed, by Claim 4 these relations are just the relations definable in terms of $\F_{|\kk}.$

In particular $\hat{\sigma}$  acts on
$\X_i(\F)$ of $T(\X(\F))$  as an automorphism of $T(\X(\F)).$ 
Now we want to
 extend the action $\hat{\sigma}$ to the whole of $\XX^{et}_ {  \{ \pp_i\}}(\F).$ Note that by (\ref{pi}) the sequence of maps $\mathbf{j_i}$  is definable in $\XX^{et}_ {  \{ \pp_i\}}(\F),$ It follows that $\U_i(\F) \subseteq \dcl (\X_i(\F))$ and thus the elementary monomorphism $\hat{\sigma}$ extends uniquely to all the sorts $\U_i(\F).$  Now the extension of  $\hat{\sigma}$ to $\U(\F)$ follows from the fact that $\U(\F)$ is the projective limit of the $\U_i(\F)$ along $\bar{\pp}_i;$
 each $u\in \U(\F)$ is the limit of the sequence $\bar{\pp}_i(u)\in \U_i(\F).$ 

Set $s(\sigma):=\hat{\sigma}.$ 
 Claim proved.
 
 \medskip
 
 It follows that $s$ is a homomorphism of $\Aut(\F_{|\kk})$ into 
$\Aut\, \XX^{et}_ {  \{ \pp_i\}}(\F)\subset \Aut\, \XX^{et}(\F). $
Thus we have $$s:  \Aut(\F_{|\kk})\to  \Aut\, \XX^{et}(\F)$$
a section associated with $T(\X).$ \footnote{Need also that the tower
$\Deck(\X_l/\X)$ has unique presentation in the form $\Gamma/\Delta_i.$}   

$\Box$
\epk
\section{Abelian varieties}
Let $\X$ be an abelian variety of dimension $g$ over $\kk,$ (in particular, $\X(\kk)\neq \emptyset$) and $\mathrm{J}(\X)$ the Jacobi variety of $\X.$

Our aim here is to construct a class of non-isomorphic cofinal towers $T(\X)$ over $\kk.$ 

\bpk 
For $n\in \N$ and $e\in \X(\kk)$
define the map $$[n]_e: \X\to \X; \ e+x\mapsto e+n\cdot x.$$  

Also fix an element $o\in \X(\kk)$ and let $$\mathcal{E}_\X=\X(\kk)^\N=\{ \{ e_i\in \X(\kk): i\in \N\}, e_0=o\}$$  the set of all sequences of elements of $\X(\kk)$ beginning with $o.$

For each $\mathbf{e}=\{ e_i\}\in \mathcal{E}_\X$ set 
$$\X_0=\X,\ \X_i=\X\mbox{ and }\mathbf{j}_{i}:=[i]_{e_i}; \X_i\to \X_{i-1}, i\in \N.$$
Clearly, $$\Deck(\X_i/\X_{i-1}), \ \Deck(\X_i/\X)\subset \mathrm{J}(\X),$$
$$\Deck(\X_i/\X_{i-1})\cong(\Z/i\Z)^{2g}, \ \ \Deck(\X_i/\X)\cong(\Z/i!\Z)^{2g},$$
the $2g$-cartesian powers of cyclic groups of orders $i$ and $i!$ respectively. 

It follows, 
$$T_\mathbf{e}(\X):\ \ \X \leftarrow_{j_{e_1}} \X_1\leftarrow_{j_{e_2}} \X_2 \leftarrow\ldots \X_i\leftarrow_{j_{e_{i+1}}} \X_{i+1}\leftarrow_{j_{e_{i+2}}} \ldots $$
is a cofinal $\Delta_i$- tower over $\kk$ for $$\Delta_i=i!\cdot \Gamma,\mbox{ where }\Gamma=\hat{\pi}_1^{top}(\X(\C)).$$
\epk
\bpk {\bf Lemma.} {\em Suppose  $T_\mathbf{e}(\X)\cong T_{\mathbf{e}'}(\X).$ Then
$(e_i-e'_i)\in \mathrm{Tors}\,\mathrm{J}(\X)$ for all $i\in \N.$  
}

{\bf Proof.} Let $f_i: \X_i\to \X'_i,$ $i\in \N,$ be the system of isomorphisms which realise the isomorphism $T_\mathbf{e}(\X)\cong T_{\mathbf{e}'}(\X).$ By definitions, 
$\X_i=\X'_i=\X$ \be \label{eq1}  f_{i-1}\circ [i]_{e_i}=[i]_{e'_i}\circ f_i\ee
Note that $f_i$ can be seen also as an isomorphism of \'etale covers
$\X_i\to \X_0$ and $\X'_i\to \X_0$ given by compositions $j_1\circ \ldots\circ j_i$ and $j'_1\circ \ldots\circ j'_i$, respectively.
It follows that $f_i$ has the form $f_i(x)= x+t_i$ for some $t\in \mathrm{Tors}(\mathrm{J}(\X)).$
 


Applying both sides of (\ref{eq1}) to $x$ we get , for $i=1,2,\ldots$
$$ f_{i-1}(i(x-e_i)+e_i)=    i\cdot(f_{i}(x)-e'_i)+e'_i.$$
in particular,
$$ f_{i-1}(e_i)=    i\cdot(f_{i}(e_i)-e'_i)+e'_i$$
and so $$e_i+t_{i-1}=i\cdot (e_i+t_i-e'_i)+e'_i$$
and finally
$$(i-1)(e_i-e'_i)=t_{i-1}-i\cdot t_i\in  \mathrm{Tors}(\mathrm{J}(\X)).$$
It follows $e_i-e'_i\in  \mathrm{Tors}(\mathrm{J}(\X)).$ $\Box$

\epk
\bpk {\bf Corollary.} {\em Assume that the group of $\kk$-rational points of $\X$ contains non-torsion points. Then there are continuum-many non-isomorphic towers $T_e(\X)$ and respectively continuum-many non-conjugated sections of the projection $\pi_1(\X)\to \Gal_\kk.$}

\epk

\end{document}